\documentclass{amsproc}
\usepackage{hyperref}

\newtheorem{theorem}{Theorem}[section]
\newtheorem{lemma}[theorem]{Lemma}
\newtheorem{proposition}[theorem]{Proposition}
\numberwithin{equation}{section}

\newcommand{\A}{\mathbb{A}}

\newcommand{\I}{\mathbb{I}}
\newcommand{\N}{\mathbb{N}}

\newcommand{\T}{\mathbb{T}}
\newcommand{\Z}{\mathbb{Z}}

\newcommand{\cB}{\mathcal{B}}
\newcommand{\cW}{\mathcal{W}}
\newcommand{\ind}{\operatorname{ind}}

\begin{document}

\title[Asymptotics of Toeplitz Determinants]
{Asymptotics of Toeplitz Determinants Generated \\
by Functions with Fourier Coefficients\\
in Weighted Orlicz Sequence Classes}

\author{Alexei Yu. Karlovich}
\address{%
Universidade do Minho,
Departamento de Matem\'atica,
Escola de Ci\^encias,
Campus de Gualtar,
4710-057, Braga, Portugal}
\email{oleksiy@math.uminho.pt}

\curraddr{%
Departamento de Matem\'atica,
Instituto Superior T\'ecnico,
Av. Rovisco Pais 1,
1049-001, Lisbon, Portugal}
\email{akarlov@math.ist.utl.pt}

\thanks{This work is supported by Centro de Matem\'atica da Universidade do Minho
(Portugal) and by the Portuguese Foundation of Science and Technology through the
research program POCTI}
\subjclass[2000]{Primary 47B35, 46E30; Secondary 15A15, 47A68}
\date{June 16, 2006}
\keywords{Toeplitz matrix, Toeplitz operator, Hankel operator,
Wiener-Hopf factorization, weighted Orlicz sequence space,
strong Szeg\H{o} limit theorem, Cauchy index}
\begin{abstract}
We prove asymptotic formulas for Toeplitz determinants generated by functions with
sequences of Fourier coefficients belonging to weighted Orlicz sequence classes.
We concentrate our attention on the case of nonvanishing generating functions with
nonzero Cauchy index.
\end{abstract}
\maketitle
\section{Introduction and main results}
\subsection{Strong Szeg\H{o}'s limit theorem for positive generating functions}
Let $\T$ be the unit circle. For a complex-valued function $a\in L^1(\T)$, let
$\{a_k\}_{k=-\infty}^\infty$ be the sequence of the Fourier coefficients of $a$,
\[
a_k:=\frac{1}{2\pi}\int_0^{2\pi}a(e^{i\theta})e^{-ik\theta}d\theta.
\]
Consider the determinants $D_n(a)$ of
the finite Toeplitz matrices $T_n(a)$,
\[
D_n(a)=\det T_n(a)=\det (a_{j-k})_{j,k=0}^n
\quad
(n\in\Z_+),
\]
where, as usual, $\Z_+:=\N\cup\{0\}$ and $\N:=\{1,2,\dots\}$.
In 1952, Gabor Szeg\H{o} \cite{S52} proved that if $a$ is a positive
function with H\"older continuous derivative, then
\begin{equation}\label{eq:Szego}
D_n(a)=G(a)^{n+1}E(a)\{1+o(1)\}
\quad\mbox{as}\quad n\to\infty,
\end{equation}
where
\[
G(a):=
e^{(\log a)_0}
=
\exp\left(\frac{1}{2\pi}\int_0^{2\pi}\log a(e^{i\theta})d\theta\right)
\]
and
\[
E(a):=
\exp\left(
\sum_{k=1}^\infty k(\log a)_k(\log a)_{-k}
\right)
\]
with $(\log a)_k$ the Fourier coefficients of $\log a$.
Basor \cite{B85} writes: ``It is interesting to note that
this formula was an important aspect of Lars Onsager's derivation
for the spontaneous magnetization of a two-dimensional Ising lattice.
The formula, for some special $a$, was proposed to Szeg\H{o} by
S.~Kakutani, who heard it from Onsager". For the importance of this
asymptotic formula in the Ising model, see \cite{MW73} and also
\cite[Section~5.2]{BS99}. The smoothness
conditions needed by Szeg\H{o} were subsequently relaxed by many authors
including Kac, Baxter, Hirschman, Krein, Devinatz, and others. Finally,
Ibragimov proved in 1968 that (\ref{eq:Szego}) is true if the constants
$E(a)$ and $G(a)$ are well defined and $a$ is positive.
For several different proofs of this beautiful result, see Simon \cite[Chap.~6]{S05}.
\subsection{Nonvanishing functions with zero Cauchy index}
Let $C(\T)$ be the Banach algebra of all complex-valued continuous functions
with the maximum norm. We will denote the \textit{Cauchy index} of a function
$a\in C(\T)$ by $\ind a$. Baxter (1963) and Hirschman (1966) were the first
to replace the positivity of $a$ in (\ref{eq:Szego}) by the condition that
\begin{equation}\label{eq:index-zero}
a(t)\ne 0 \quad\mbox{for all}\quad t\in\T\quad\mbox{and}\quad \ind a=0.
\end{equation}
To formulate results of this kind precisely, we will define some smoothness
classes in terms of the decay of the Fourier coefficients.
\subsection{The Wiener algebra}
Let $W$ be the \textit{Wiener algebra} of all complex-valued functions $a$ on
$\T$ of the form
\[
a(t)=\sum_{k=-\infty}^\infty a_k t^k
\quad (t\in\T)
\quad\mbox
{for which}\quad
\|a\|_W:=\sum_{k=-\infty}^\infty |a_k|<\infty.
\]
It is well known that $W$ is a Banach algebra under the norm $\|\cdot\|_W$
and that $W$ is continuously embedded into $C(\T)$.
\subsection{Functions with Fourier coefficients in weighted Orlicz sequence classes}
\label{sec:Orlicz-class}
Let $p:[0,\infty)\to[0,\infty)$ be a right-continuous non-decreasing function
such that $p(0)=0$, $p(t)>0$ for $t>0$, and $\lim\limits_{t\to\infty} p(t)=\infty$.
Then the function $q(s)=\sup\{t : p(t)\le s\}$ (defined for $s\ge 0$) has the same
properties as the function $p$. The convex functions $\Phi$ and $\Psi$ defined by
the equalities
\[
\Phi(x):=\int_0^x p(t) dt,
\quad
\Psi(x):=\int_0^x q(s)ds
\quad (x\ge 0)
\]
are called \textit{complementary $N$-functions} (see, e.g.,
\cite[Section~1.3]{KR61},
\cite[Ch.~8]{M89},
\cite[Section~13]{M83}).
An $N$-function $\Phi$ is said to satisfy the $\Delta_2^0$-condition if
\[
\limsup_{x\to 0}\frac{\Phi(2x)}{\Phi(x)}<\infty.
\]

Any sequence $\{\nu_k\}_{k=0}^\infty$ of positive numbers is called a weight
sequence. We denote by $\cW$ the set of all weight sequences
$\{\nu_k\}_{k=0}^\infty$  such that
\begin{enumerate}
\item[(a)]
$\nu_0=1$;
\item[(b)]
$\nu_{k-1}\le\nu_k$ for $k\in\N$;
\item[(c)]
$\{\nu_k\}_{k=0}^\infty$ satisfies
the $\Delta_2^\N$-condition, that is, there is a constant $C_\nu\in(0,\infty)$
such that $\nu_{2k}\le C_\nu\nu_k$ for $k\in\N$.
\end{enumerate}
It is easy to see that $C_\nu\ge 1$.

Given $N$-functions $\Phi,\Psi$ and weight sequences $\varphi=\{\varphi_k\}_{k=0}^\infty$,
$\psi=\{\psi_k\}_{k=0}^\infty$, denote by $F\ell_{\varphi,\psi}^{\Phi,\Psi}$
the set of all complex-valued functions $a\in L^1(\T)$ satisfying
\begin{equation}\label{eq:Orlicz-class}
\sum_{k=1}^\infty \Phi(|a_{-k}|\varphi_k)
+
\sum_{k=0}^\infty \Psi(|a_k|\psi_k)<\infty.
\end{equation}
If $\Phi(t)=t^p/p$, $\Psi(t)=t^q/q$ with $p,q\ge 1$ or $\varphi(k)=(k+1)^\alpha$,
$\psi(k)=(k+1)^\beta$ with $\alpha,\beta\ge 0$, we will simply write
$F\ell_{\varphi,\psi}^{p,q}$ or $F\ell_{\alpha,\beta}^{\Phi,\Psi}$ instead
of $F\ell_{\varphi,\psi}^{\Phi,\Psi}$.
\subsection{The strong Szeg\H{o} limit theorem \`a la Hirschman}
Hirschman \cite{H66} proved that $W\cap F\ell_{1/2,1/2}^{2,2}$ is a Banach
algebra under pointwise multiplication and that if $a\in W\cap F\ell_{1/2,1/2}^{2,2}$
is a complex-valued function satisfying (\ref{eq:index-zero}),
then (\ref{eq:Szego}) is valid. B\"ottcher and Silbermann \cite{BS80}
(see also \cite[Corollary~10.45]{BS06}) proved a
non-symmetric version of the latter result. Suppose $\alpha\in[0,1]$, $p>1$,
and, in addition, $W\cap F\ell_{\alpha,1-\alpha}^{p,p/(p-1)}$ is a Banach algebra
under pointwise multiplication. If $a\in W\cap F\ell_{\alpha,1-\alpha}^{p,p/(p-1)}$
satisfies (\ref{eq:index-zero}), then (\ref{eq:Szego}) holds.
They conjectured also that $W\cap F\ell_{\alpha,\beta}^{p,q}$
is always a Banach algebra whenever $\alpha\ge 0$, $\beta\ge 0$ and $p\ge 1$, $q\ge 1$.
This conjecture was proved by Horbach in 1984 (see \cite[Theorem~6.54]{BS06}).

Using the ideas of the proof of \cite[Theorem~6.54]{BS06}, the author \cite{K04} proved
that if $\Phi,\Psi$ are $N$-functions both satisfying the $\Delta_2^0$-condition
and $\varphi,\psi\in\cW$, then $W\cap F\ell_{\varphi,\psi}^{\Phi,\Psi}$ is a Banach
algebra under pointwise multiplication and an appropriate norm (see also
Lemma~\ref{le:algebra}). Using this result and
ideas of the proof of \cite[Corollary~10.45]{BS06}, the author and Santos \cite{KS05}
(see also \cite[Section~10.46]{BS06}) obtained the following version of the
strong Szeg\H{o} limit theorem.
\begin{theorem}\label{th:KS}
Suppose $\Phi,\Psi$ are complementary $N$-functions both satisfying
the $\Delta_2^0$-condition, $\varphi=\{\varphi_k\}_{k=0}^\infty$,
$\psi=\{\psi_k\}_{k=0}^\infty$ are weight sequences in $\cW$,
and there exists a constant $M\in(0,\infty)$ such that
$k\le M\varphi_k\psi_k$ for all $k\in\Z_+$.
If $a\in W\cap F\ell_{\varphi,\psi}^{\Phi,\Psi}$ satisfies
{\rm (\ref{eq:index-zero})}, then {\rm (\ref{eq:Szego})} holds.
\end{theorem}
\subsection{Further results on asymptotics of Toeplitz determinants}
There exist generalizations of the strong Szeg\H{o} limit theorem into different
directions. For instance, Widom \cite{W76} extended (\ref{eq:Szego}) to the case of
matrix-valued generating functions and B\"ottcher and Silbermann \cite{BS94}
extended further this result to the case of operator-valued generating functions.

On the other hand, if one of the assumptions in (\ref{eq:index-zero}) is not satisfied,
then the asymptotic formula (\ref{eq:Szego}) may fail. The asymptotics of Toeplitz
determinants with generating functions which do not satisfy one of the conditions
in (\ref{eq:index-zero}) were first considered by Fisher and Hartwig \cite{FH68,FH69}.
In particular, they conjectured the asymptotic behavior of  Toeplitz determinants with
some interesting singular generating functions (the so-called Fisher-Hartwig generating
functions). Later their conjecture was proved in many particular cases in papers
by Basor, B\"ottcher, Silbermann, Widom, and others. Basor and Tracy \cite{BT91} provided
a series of counterexamples to the Fisher-Hartwig conjecture and stated a generalized
conjecture. Ehrhardt \cite{E97} made a significant progress in proving the Basor-Tracy
conjecture. We refer for the status of this problem to Ehrhardt's survey \cite{E01} and
also to \cite[Sections 10.57--10.80]{BS06}.

Many aspects of asymptotic behavior of Toeplitz determinants are considered
in the monographs by Grenander and Szeg\H{o} \cite{GS58}, B\"ottcher and Grudsky
\cite{BG05}, B\"ottcher and Silbermann \cite{BS83,BS99,BS06}, and Simon \cite{S05},
where the reader can find further results, historical remarks, and references.
\subsection{The case of a nonvanishing Cauchy index: known results}
In this paper we embark on the case of nonsingular generating functions with
nonzero Cauchy index. Fisher and Hartwig \cite{FH68,FH69} proved asymptotic formulas
for Toeplitz determinants generated by functions of the form $a(t)t^\kappa$, where
$a$ satisfies (\ref{eq:index-zero}) and $\kappa\in\Z$. Clearly, in that case
$\ind[a(t)t^\kappa]=\kappa$. Their results were extended to the case of matrix-valued
generating functions by B\"ottcher and Silbermann \cite[Theorem~14]{BS80}.
They considered generating functions from the (matrix version of)
H\"older-Zygmund spaces. Let $\omega_2(g,s)$ denote the modulus of smoothness of order
$2$ of a bounded function $g$ on $\T$. Given $\gamma>0$, write $\gamma=m+\delta$, where
$m\in\Z_+$ and $\delta\in(0,1]$. The H\"older-Zygmund space $C^\gamma=C^\gamma(\T)$ is
defined (see, e.g., \cite[Section~3.5.4]{ST87}) by
\[
C^\gamma:=\left\{
f\in C(\T): f^{(j)}\in C(\T),\ 1\le j\le m, \
\sup_{s>0}\frac{\omega_2(f^{(m)},s)}{s^\delta}<\infty
\right\},
\]
where $f^{(j)}$ is the derivative of order $j$ of $f$.
Consider the Hardy spaces
\begin{eqnarray*}
H_\pm^\infty
&:=&
\big\{a\in L^\infty(\T) : a_{\mp n}=0
\mbox{ for } n\in\N\big\}.
\end{eqnarray*}
\begin{theorem}\label{th:nonvanishing-1}
Let $\gamma>0$ and $\A=C^\gamma$. If $a\in\A$ satisfies {\rm(\ref{eq:index-zero})},
then

\noindent
{\rm (a)}
there exist functions $a_-$ and $a_+$ such that $a=a_-a_+$ and
\[
a_-^{\pm 1}\in \A\cap H_-^\infty,
\quad
a_+^{\pm 1}\in\A\cap H_+^\infty;
\]

\noindent
{\rm (b)} the following refinement of {\rm(\ref{eq:Szego})} is valid:
\[
D_n(a)=G(a)^{n+1}E(a)\{1+\delta_1(n)\} \quad\mbox{as}\quad n\to\infty
\]
with $\delta_1(n)=O(n^{1-2\gamma})$;

\medskip
\noindent
{\rm (c)}
if we put $b=a_-a_+^{-1}$ and $c=a_-^{-1}a_+$, then for every $\kappa\in\N$,
\begin{eqnarray*}
&&
D_n[a(t)t^{-\kappa}]
=
G(a)^{n+1} E(a) \{1+\delta_1(n)\}
\\
&&
\quad\quad
\times(-1)^{(n+\kappa)\kappa}G(c)^\kappa
\left\{\det\left(\begin{array}{ccc}
b_{n+1} &\dots & b_{n-\kappa+2}\\
\vdots & \ddots &\vdots\\
b_{n+\kappa} &\dots & b_{n+1}
\end{array}\right)+\delta_2(n)
\right\}
\end{eqnarray*}
and
\begin{eqnarray*}
&&
D_n[a(t)t^{\kappa}]
=
G(a)^{n+1} E(a)\{1+\delta_1(n)\}
\\
&&
\quad\quad
\times
(-1)^{(n+\kappa)\kappa}G(b)^\kappa
\left\{\det\left(\begin{array}{ccc}
c_{-n-1} &\dots & c_{-n+\kappa-2}\\
\vdots & \ddots &\vdots\\
c_{-n-\kappa} &\dots & c_{-n-1}
\end{array}\right)+\delta_3(n)
\right\}
\end{eqnarray*}
as $n\to\infty$ with $\delta_1(n)=O(n^{1-2\gamma})$
and $\delta_2(n)=\delta_3(n)=O(n^{-3\gamma})$.
\end{theorem}

The representation in (a) is called a \textit{Wiener-Hopf factorization} of
$a$ in the algebra $\A$.

The proof of part (a) can be found in \cite[Section 6.25(vi)]{PS91}.
Part (b) follows from \cite[Sections~6.18--6.20]{BS83} or
\cite[Theorems~10.35 and 10.37]{BS06}. Part (c) follows from part (a)
and from \cite[Theorem~6.24]{BS83} or \cite[Theorem~10.47]{BS06}.

An alternative approach to the asymptotics of Toeplitz determinants
with scalar-valued generating functions of nonvanishing Cauchy index is suggested in a series
of papers by Carey and Pincus \cite{CP99,CP01,CP06}. In \cite[Theorem~A]{CP06},
they found an exact formula for Toeplitz determinants generated by functions
of nonvanishing Cauchy index. Their approach is based on a heavy use of results
and methods of algebraic $K$-theory of the algebras of operators having trace class
commutators and it is by no means elementary. Very recently B\"ottcher and Widom
\cite{BW06} have found an elementary proof of the above mentioned exact formula
and, therefore, have obtained a new proof of Theorem~\ref{th:nonvanishing-1}(c).

The author \cite{K06c} (see also \cite{K06a,K06b}) has obtained an analog of
Theorem~\ref{th:nonvanishing-1} for matrix-valued generating functions in
weighted Wiener algebras. Its scalar version for weight sequences in $\cW$
reads as follows.
\begin{theorem}\label{th:nonvanishing-2}
Suppose $\varphi=\{\varphi_j\}_{j=0}^\infty$, $\psi=\{\psi_j\}_{j=0}^\infty$
are weight sequences in $\cW$ and $\sum_{j=1}^\infty[\varphi_j\psi_j]^{-1}<\infty$.
If $a\in W\cap F\ell_{\varphi,\psi}^{1,1}$ satisfies {\rm(\ref{eq:index-zero})},
then the results of Theorem~{\rm\ref{th:nonvanishing-1}(a)--(c)} are valid with
$\A=W\cap F\ell_{\varphi,\psi}^{1,1}$ and
\begin{eqnarray*}
\delta_1(n)&=&o\left(\sum_{j=n+1}^\infty[\varphi_{j+1}\psi_{j+1}]^{-1}\right),
\\
\delta_2(n)&=&o\big(\varphi_{n+\kappa+1}^{-1}\psi_{n+\kappa+1}^{-2}\big),
\\
\delta_3(n)&=&o\big(\varphi_{n+\kappa+1}^{-2}\psi_{n+\kappa+1}^{-1}\big).
\end{eqnarray*}
\end{theorem}

This result follows from \cite[Theorem~1.2]{K06c} and \cite[Proposition~22]{K06a}.

In particular, if $\A=W\cap F\ell_{\alpha,\beta}^{1,1}$ with $\alpha,\beta\ge 0$ and
$\alpha+\beta>1$, then the results of Theorem~{\rm\ref{th:nonvanishing-1}(a)--(c)}
are valid with
\[
\delta_1(n)=o(n^{1-\alpha-\beta}),
\quad
\delta_2(n)=o(n^{-\alpha-2\beta}),
\quad
\delta_3(n)=o(n^{-2\alpha-\beta}).
\]
\subsection{The case of a nonvanishing Cauchy index: new results}
Our main result is the following version of Theorem~\ref{th:nonvanishing-1}
for generating functions in $W\cap F\ell_{\varphi,\psi}^{\Phi,\Psi}$.
\begin{theorem}\label{th:main}
Suppose $\Phi,\Psi$ are complementary $N$-functions both satisfying
the $\Delta_2^0$-condition, $\varphi=\{\varphi_k\}_{k=0}^\infty$,
$\psi=\{\psi_k\}_{k=0}^\infty$ are weight sequences in $\cW$,
and there exists a constant $M\in(0,\infty)$ such that
$k\le M\varphi_k\psi_k$ for all $k\in\Z_+$.
If $a\in W\cap F\ell_{\varphi,\psi}^{\Phi,\Psi}$ satisfies {\rm (\ref{eq:index-zero})},
then the results of Theorem~{\rm\ref{th:nonvanishing-1}(a), (c)} are valid with
$\A=W\cap F\ell_{\varphi,\psi}^{\Phi,\Psi}$ and
\[
\delta_1(n)=o(1),
\quad
\delta_2(n)=o(1/\psi_n),
\quad
\delta_3(n)=o(1/\varphi_n).
\]
\end{theorem}
Part (a) was proved in \cite[Corollary~2.3]{K04} (see also Lemma~\ref{le:algebra}(b)).
Part (c) will be proved in Section~\ref{subsec:proof}.

If $\Phi$ and $\Psi$ are arbitrary $N$-functions, then from the well known
conditions for the embeddings of Orlicz classes (see, e.g., \cite[Theorem~3.4(c)]{M89}
or \cite[Section~8]{M83} it follows that
$F\ell_{\varphi,\psi}^{1,1}\subset F\ell_{\varphi,\psi}^{\Phi,\Psi}$.
Thus the class $W\cap F\ell_{\varphi,\psi}^{\Phi,\Psi}$
of generating functions in Theorem~\ref{th:main} is larger than the class
$W\cap F\ell_{\varphi,\psi}^{1,1}$ of generating functions in
Theorem~\ref{th:nonvanishing-2}. On the other hand, one has better speed of
convergence in Theorem~\ref{th:nonvanishing-2} than in Theorem~\ref{th:main}.

The paper is organized as follows. In Section~\ref{sec:aux} we give definitions
of weighted Orlicz sequence spaces and classes, define the norm in the class
$F\ell_{\varphi,\psi}^{\Phi,\Psi}$ and collect necessary results about
this class and its intersection with the Wiener algebra $W$.
Finally we formulate basic results on Toeplitz and Hankel operators acting on
non-weighted Orlicz sequence spaces. In Section~\ref{sec:proof}
we give the proof of Theorem~\ref{th:main}(c) following the approach
developed by B\"ottcher and Silbermann in \cite[Section~8]{BS80}
(see also \cite[Section~10.47]{BS06}). We apply Jacobi's
theorem on the first step. It allows us to represent $D_{n-\kappa}[a(t)t^\kappa]$
as the product of $D_n(a)$ (it is treated by Theorem~\ref{th:KS}) and a minor of
$T_n^{-1}(a)$. This minor is represented in Lemma~\ref{le:BS2} as a sum
of two terms. The first term  gives the leading term in the asymptotic formula
after some computations. If the generating function is sufficiently smooth
(in our case $a\in W\cap F\ell_{\varphi,\psi}^{\Phi,\Psi}$), then the norm of the
second term is asymptotically small (in our case $o(1/\varphi_{n-\kappa})$
or $o(1/\psi_{n-\kappa})$). Gathering all these pieces together, we finish the
proof in Section~\ref{subsec:proof}.
\section{Auxiliary results}\label{sec:aux}
\subsection{One result from numerical linear algebra}
The following result was stated in \cite[Section~10.47]{BS06} without a proof,
a sketch of the proof can be found in \cite[Proposition~2.3]{K06c}.
\begin{proposition}\label{pr:LA}
Suppose $\{\gamma_n\}_{n=1}^\infty$ is a sequence of positive numbers
and $\{A_n\}_{n=1}^\infty$, $\{B_n\}_{n=1}^\infty$ are sequences of $m\times m$
matrices. If $\sup\limits_{n\in\N}\|A_n\|<\infty$
and $\|A_n-B_n\|=o(\gamma_n)$ as $n\to\infty$, where $\|\cdot\|$ is any matrix norm, then
\[
\det A_n=\det B_n+o(\gamma_n)\quad\mbox{as}\quad n\to\infty.
\]
\end{proposition}
\subsection{Weighted Orlicz sequence spaces and classes}
Let $\I$ be either $\N$ or $\Z_+$. Suppose $\Phi$
is an $N$-function and $\varphi=\{\varphi_k\}_{k=0}^\infty$ is a weight sequence.
The set $\ell^\Phi_\varphi(\I)$ of all sequences $c=\{c_k\}_{k\in\I}$
of complex numbers such that
\begin{equation}\label{eq:Orlicz-def}
\sum_{k\in\I}\Phi\left(\frac{|c_k|\varphi_k}{\lambda}\right)<\infty
\end{equation}
for some $\lambda=\lambda(c)>0$ is a Banach space when equipped with the norm
\[
\|c\|_{\ell^\Phi_\varphi(\I)}
=
\inf\left\{\lambda>0\ :\
\sum_{k\in\I}\Phi\left(\frac{|c_k|\varphi_k}{\lambda}\right)\le 1
\right\}.
\]
The space $\ell^\Phi_\varphi(\I)$ is called a \textit{weighted Orlicz
sequence space}.  If $\varphi_k=1$ for all $k\in\Z_+$, then we will
simply write $\ell^\Phi(\I)$ instead of $\ell^\Phi_\varphi(\I)$ and we
will say that $\ell^\Phi(\I)$ is an \textit{Orlicz sequence space}.
The set $\widetilde{\ell}^\Phi_\varphi(\I)$
of all sequences $c=\{c_k\}_{k\in\I}$ of complex numbers for which
(\ref{eq:Orlicz-def}) holds with $\lambda=1$ is called a
\textit{weighted Orlicz sequence class}.
Weighted Orlicz sequence spaces are the partial case of so-called
Musie\-lak-Orlicz sequence spaces (= modular sequence spaces).
Good sources for the theory of Musielak-Orlicz sequence spaces
are \cite[Section~4.d]{LT77} and \cite{M83}; for Orlicz
sequence spaces see also \cite{M89}; for Orlicz function spaces
over finite measure spaces, see \cite{KR61}.

Applying the results of \cite[Theorem~8.14(b)]{M83}
(see also \cite[Proposition~4.d.3]{LT77}) to the sequence of
$N$-functions $\Phi_k(x)=\Phi(x\varphi_k)$, we get the following.
\begin{theorem}\label{th:Orlicz-class}
Suppose $\Phi$ is an $N$-function satisfying the $\Delta_2^0$-condition and
$\varphi=\{\varphi_k\}_{k=0}^\infty$ is a weight sequence. Then
$\ell^\Phi_\varphi(\I)=\widetilde{\ell}^\Phi_\varphi(\I)$.
\end{theorem}
\subsection{The algebra $W\cap F\ell_{\varphi,\psi}^{\Phi,\Psi}$}
Let $\Phi,\Psi$ be $N$-functions and let
$\varphi=\{\varphi_k\}_{k=0}^\infty$, $\psi=\{\psi_k\}_{k=0}^\infty$
be weight sequences. If $\Phi$ and $\Psi$ both satisfy the
$\Delta_2^0$-condition, then from Theorem~\ref{th:Orlicz-class} it follows
that the set $F\ell^{\Phi,\Psi}_{\varphi,\psi}$ of all
functions $a\in L^1(\T)$ with the Fourier coefficients $\{a_k\}_{k=-\infty}^\infty$
satisfying (\ref{eq:Orlicz-class}) is a Banach space with respect to the norm
\[
\|a\|_{F\ell^{\Phi,\Psi}_{\varphi,\psi}}
:=
\|\{a_{-k}\}_{k\in\N}\|_{\ell^\Phi_\varphi(\N)}
+
\|\{a_k\}_{k\in\Z_+}\|_{\ell^\Psi_\psi(\Z_+)}.
\]
For $a\in L^1(\T)$ and $n\in\N$, put
\[
a^{(n)}(t)=\sum_{k=-n}^n a_kt^k\quad (t\in\T).
\]
\begin{lemma}\label{le:properties}
{\rm(see \cite[Lemma~3, Proposition~2]{KS05}).}
Let $\Phi,\Psi$ be $N$-functions both satisfying the $\Delta_2^0$-condition and let
$\varphi=\{\varphi_k\}_{k=0}^\infty$, $\psi=\{\psi_k\}_{k=0}^\infty$
be weight sequences.

{\rm (a)}
There is a constant $C(\varphi,\psi,\Phi,\Psi)>0$
depending only on $\varphi,\psi$ and $\Phi,\Psi$ such that
for all $a\in F\ell_{\varphi,\psi}^{\Phi,\Psi}$,
\[
\|\overline{a}\|_{F\ell_{\psi,\varphi}^{\Psi,\Phi}}
\le
C(\varphi,\psi,\Phi,\Psi)
\|a\|_{F\ell_{\varphi,\psi}^{\Phi,\Psi}}.
\]

{\rm (b)}
If $a\in F\ell_{\varphi,\psi}^{\Phi,\Psi}$, then
\[
\lim_{n\to\infty}\|a-a^{(n)}\|_{F\ell_{\varphi,\psi}^{\Phi,\Psi}}=0.
\]
\end{lemma}
We equip the set $W\cap F\ell_{\varphi,\psi}^{\Phi,\Psi}$ with the norm
\begin{equation}\label{eq:norm}
\|a\|_{W\cap F\ell_{\varphi,\psi}^{\Phi,\Psi}}:=
\|a\|_W+\|a\|_{F\ell_{\varphi,\psi}^{\Phi,\Psi}}.
\end{equation}
The following result generalizes Horbach's theorem \cite[Theorem~6.54]{BS06}.
It was recently proved in \cite{K04} in a slightly more general form.
\begin{lemma}\label{le:algebra}
Let $\Phi,\Psi$ be $N$-functions both satisfying the $\Delta_2^0$-condition and let
$\varphi=\{\varphi_k\}_{k=0}^\infty$, $\psi=\{\psi_k\}_{k=0}^\infty$
be weight sequences in $\cW$.

{\rm (a)}
If $a,b\in W\cap F\ell_{\varphi,\psi}^{\Phi,\Psi}$, then
\[
\|ab\|_{W\cap F\ell_{\varphi,\psi}^{\Phi,\Psi}}
\le
(1+2C_\varphi+2C_\psi)
\|a\|_{W\cap F\ell_{\varphi,\psi}^{\Phi,\Psi}}
\|b\|_{W\cap F\ell_{\varphi,\psi}^{\Phi,\Psi}}.
\]

{\rm (b)}
If $a\in W\cap F\ell_{\varphi,\psi}^{\Phi,\Psi}$ satisfies {\rm(\ref{eq:index-zero})},
then $a$ has a logarithm in $W\cap F\ell_{\varphi,\psi}^{\Phi,\Psi}$. If we let
for $t\in\T$,
\[
a_-(t):=\exp\left(\sum_{k=1}^\infty (\log a)_{-k}t^{-k}\right),
\quad
a_+(t):=\exp\left(\sum_{k=0}^\infty(\log a)_kt^k\right),
\]
then $a=a_-a_+$ and
\[
a_-^{\pm 1}\in (W\cap F\ell_{\varphi,\psi}^{\Phi,\Psi})\cap H_-^\infty,
\quad
a_+^{\pm 1}\in (W\cap F\ell_{\varphi,\psi}^{\Phi,\Psi})\cap H_+^\infty.
\]
\end{lemma}
\subsection{Toeplitz and Hankel operators on Orlicz sequence spaces}
Let $a$ be a function in $L^1(\T)$ with the Fourier coefficients
$\{a_k\}_{k=-\infty}^\infty$ and let $\{c_k\}_{k=-\infty}^\infty$
be a sequence of complex numbers. We formally define the \textit{Laurent operator}
with the symbol $a$ by
\[
L(a):\{c_k\}_{k=-\infty}^\infty\mapsto
\Big\{\sum_{k=-\infty}^\infty a_{j-k}c_k\Big\}_{j=-\infty}^\infty
\]
and the operators $P,Q$, and $J$ as follows
\[
(Pc)_k:=\left\{\begin{array}{ccc}
c_k & \mbox{for} & k\ge 0,\\
0   & \mbox{for} & k<0,
\end{array}\right.
\quad
(Qc)_k:=\left\{\begin{array}{ccc}
0 & \mbox{for} & k\ge 0,\\
c_k   & \mbox{for} & k<0,
\end{array}\right.
\quad (Jc)_k=c_{-k-1}.
\]
For $t\in\T$, put $\widetilde{a}(t):=a(1/t)$. Define \textit{Toeplitz operators}
\[
T(a):=PL(a)P|\operatorname{Im}P,
\quad
T(\widetilde{a}):=JQL(a)QJ|\operatorname{Im}P
\]
and \textit{Hankel operators}
\[
H(a):=PL(a)QJ|\operatorname{Im}P,
\quad
H(\widetilde{a}):=JQL(a)P|\operatorname{Im}P.
\]

For the spaces $\ell^p(\Z_+)$, the following result is well known
(see \cite[Chap.~2]{BS83}, \cite[Chap.~1]{BS99}, \cite[Chap.~2]{BS06}),
for Orlicz sequence spaces the proofs are actually the same
(see \cite[Section~3]{KS05}).
\begin{lemma}\label{le:basic-Toeplitz}
Suppose $\Psi$ is an $N$-function.

{\rm (a)} If $a\in W$, then the operators $T(a)$ and $T(\widetilde{a})$
are bounded on $\ell^\Psi(\Z_+)$ and the operators $H(a)$ and $H(\widetilde{a})$
are compact on $\ell^\Psi(\Z_+)$.

{\rm (b)} If $a,b\in W$, then
\begin{eqnarray*}
T(ab)&=&T(a)T(b)+H(a)H(\widetilde{b}),
\\
H(ab)&=&T(a)H(b)+H(a)T(\widetilde{b}).
\end{eqnarray*}

{\rm (c)} If $a_-\in W\cap H_-^\infty$, then $H(a_-)=0$. If $a_+\in W\cap H_+^\infty$,
then $H(\widetilde{a_+})=0$.
\end{lemma}
\section{Proof of the main result}
\label{sec:proof}
\subsection{Application of Jacobi's theorem}
Let $\Psi$ be an $N$-function. Denote by $\cB(\ell^\Psi(\Z_+))$
the Banach algebra of all bounded linear operators on the Orlicz sequence
space $\ell^\Psi(\Z_+)$. For $n\in\Z_+$, define the operators $P_n$ and $Q_n$ by
\[
P_n:\{c_k\}_{k=0}^\infty\mapsto\{c_0,c_1,\dots,c_n,0,0,\dots\},
\quad
Q_n:=I-P_n.
\]
Obviously, $P_n,Q_n\in\cB(\ell^\Psi(\Z_+))$ and $P_n^2=P_n$, $Q_n^2=Q_n$.

We will identify the operator $P_nT(a)P_n:P_n\ell^\Psi(\Z_+)\to P_n\ell^\Psi(\Z_+)$
with the finite Toeplitz matrix $T_n(a)=[a_{j-k}]_{j,k=0}^n$, if this operator
is invertible we will simply write $T_n^{-1}(a)$ instead of $(P_nT(a)P_n)^{-1}P_n$.

Fisher and Hartwig \cite{FH68,FH69} recognized that the following result
is the key to treating the asymptotics of Toeplitz determinants generating
by functions of nonvanishing Cauchy index.
\begin{lemma}\label{le:Jacobi}
Let $\kappa\in\N$ and $n\ge\kappa$. If $a\in W$ is
such that $T_n(a)$ is invertible, then
\begin{eqnarray*}
\det\big[(P_n-P_{n-\kappa})T_n^{-1}(a)P_{\kappa-1}\big]
&=&
(-1)^{n\kappa}\frac{D_{n-\kappa}[a(t)t^{-\kappa}]}{D_n(a)},
\\
\det\big[P_{\kappa-1}T_n^{-1}(a)(P_n-P_{n-\kappa})\big]
&=&
(-1)^{n\kappa}\frac{D_{n-\kappa}[a(t)t^\kappa]}{D_n(a)}.
\end{eqnarray*}
\end{lemma}
\begin{proof}
This lemma follows from Jacobi's theorem on the conjugate minors of the
adjugate matrix (see, e.g., \cite[Chap.~I, Section 4]{G59}).
This theorem is applied to $T_n(a)$ and the $\kappa\times\kappa$ minor
standing at the left lower (resp. right upper) corner of $T_n^{-1}(a)$.
\end{proof}
\subsection{The B\"ottcher-Silbermann asymptotic analysis}
Our starting point is the following simple and important fact.
\begin{proposition}\label{pr:strong}
{\rm(see \cite[Proposition~11]{KS05}).}
If $\Psi$ is an $N$-function satisfying the $\Delta_2^0$-condition, then
the sequence $P_n$ converges strongly to the identity operator $I$ on the space
$\ell^\Psi(\Z_+)$.
\end{proposition}

In \cite[Lemma~6(a)]{KS05} we proved that the so-called \textit{finite section method}
is applicable to Toeplitz operators on Orlicz sequence spaces. One of the
equivalent forms of this property can be stated as follows (see, e.g.,
\cite[Proposition~7.3]{BS06}).
\begin{lemma}\label{le:stability}
Let $\Psi$ be an $N$-function satisfying the $\Delta_2^0$-condition.
Suppose $a_-^{\pm 1}\in W\cap H_-^\infty$, $a_+^{\pm 1}\in W\cap H_+^\infty$, and
put $a=a_-a_+$. Then the sequence $\{T_n(a)\}_{n=0}^\infty$ is stable in $\ell^\Psi(\Z_+)$,
that is, for all sufficiently large $n$, say $n\ge n_0$, the matrices $T_n(a)$
are invertible and
\begin{equation}\label{eq:BS1-1}
\sup_{n\ge n_0}\|T_n^{-1}(a)\|_{\cB(\ell^\Psi(\Z_+))}<\infty.
\end{equation}
\end{lemma}

B\"ottcher and Silbermann \cite{BS80} developed further Widom's ideas \cite{W76}
and suggested an approach to study of asymptotics of Toeplitz determinants
based on the Wiener-Hopf factorization. We formulate their key
identities in the setting of Orlicz sequence spaces (originally stated in the
setting of $\ell^2(\Z_+)$).
\begin{lemma}\label{le:BS1}
Let $\Psi$ be an $N$-function satisfying the $\Delta_2^0$-condition.
Suppose $a_-^{\pm 1}\in W\cap H_-^\infty$, $a_+^{\pm 1}\in W\cap H_+^\infty$, and
put $a=a_-a_+$, $b=a_-a_+^{-1}$, and $c=a_-^{-1}a_+$. Then for all sufficiently
large $n$, say $n\ge n_0$, the matrices $T_n(a)$ are invertible and
\begin{equation}\label{eq:BS1-2}
T_n^{-1}(a)=P_nT(a_+^{-1})P_n\left\{I-\sum_{m=0}^\infty F_{n,m}\right\}P_nT(a_-^{-1})P_n,
\end{equation}
where, for $m,n\in\N$,
\begin{equation}\label{eq:BS1-3}
F_{n,0}:=P_nT(c)Q_nT(b)P_n,
\quad
F_{n,m}:=P_nT(c)(Q_nH(b)H(\widetilde{c})Q_n)^mT(b)P_n,
\end{equation}
and the convergence in {\rm(\ref{eq:BS1-2})} is understood in the sense of
$\cB(\ell^\Psi(\Z_+))$.
\end{lemma}
\begin{proof}
This statement can be found in
\cite[Section~6.15]{BS83} or
\cite[Section~10.34]{BS06} in the case of $\ell^2(\Z_+)$.
We refer also to \cite[Lemma~8]{KS05}, where a similar formula for $P_0T_n^{-1}(a)P_0$
is proved in detail for $\ell^\Psi(\Z_+)$.

Below we give a sketch of the proof.
By Lemma~\ref{le:properties}, the operators $T(a)$ and $T(a^{-1})$ are invertible
on $\ell^\Psi(\Z_+)$ and
\begin{equation}\label{eq:BS1-4}
T^{-1}(a)=T(a_+^{-1})T(a_-^{-1}),
\quad
T^{-1}(a^{-1})=T(a_+)T(a_-).
\end{equation}
By Lemma~\ref{le:stability}, there exists a number $n_0$ such that the matrices
$T_n(a)$ are invertible for all $n\ge n_0$ and (\ref{eq:BS1-1}) is fulfilled.

Applying \cite[Proposition~7.15]{BS06}, (\ref{eq:BS1-4}) and
Lemma~\ref{le:basic-Toeplitz}(c), we obtain that the operators $Q_nT^{-1}(a)Q_n$
are invertible on $Q_n\ell^\Psi(\Z_+)$ and
\begin{eqnarray}\label{eq:BS1-5}
T_n^{-1}(a)&=&
P_nT(a_+^{-1})P_n
\\
\nonumber
&&
\times\{I-P_nT(a_-^{-1})Q_n(Q_nT^{-1}(a)Q_n)^{-1}Q_nT(a_+^{-1})P_n\}
\\
\nonumber
&&
\times P_nT(a_-^{-1})P_n.
\end{eqnarray}
By Lemma~\ref{le:basic-Toeplitz}, $T^{-1}(a)=T(a^{-1})-K$, where the operator
$K:=H(a_+^{-1})H(\widetilde{a_-^{-1}})$ is compact on $\ell^\Psi(\Z_+)$.
From the identities
\begin{equation}\label{eq:BS1-6}
P_nT(a_+^{\pm 1})Q_n=Q_nT(a_-^{\pm 1})P_n=0
\end{equation}
and
Lemma~\ref{le:basic-Toeplitz} it follows that the operators
$A_n:=Q_nT(a^{-1})Q_n|\operatorname{Im}Q_n$ are invertible and
\begin{equation}\label{eq:BS1-7}
A_n^{-1}Q_n=Q_nT(a_+)Q_nT(a_-)Q_n
\end{equation}
are uniformly bounded. Since $K$ is compact
on $\ell^\Psi(\Z_+)$, taking into account Proposition~\ref{pr:strong},
we get $\|K_n\|_{\cB(\ell^\Psi(\Z_+))}\to 0$ as $n\to\infty$, where
$K_n:=Q_nKQ_n$. Hence $\|A_n^{-1}K_n\|_{\cB(\ell^\Psi(\Z_+))}\to 0$ as $n\to\infty$
and
\begin{eqnarray*}
(Q_nT^{-1}(a)Q_n)^{-1}Q_n
&=&
(Q_nT(a^{-1})Q_n-Q_nKQ_n)^{-1}Q_n=(A_n-K_n)^{-1}Q_n
\\
&=&
(I-A_n^{-1}K_n)^{-1}A_n^{-1}Q_n
=\sum_{m=0}^\infty(A_n^{-1}K_n)^m A_n^{-1}Q_n.
\end{eqnarray*}
Combining this identity with (\ref{eq:BS1-5}), we arrive at (\ref{eq:BS1-2})
with
\begin{eqnarray*}
F_{n,0} &:=& P_nT(a_-^{-1})Q_nA_n^{-1}Q_nT(a_+^{-1})P_n,
\\
F_{n,m} &:=& P_nT(a_-^{-1})Q_n(A_n^{-1}K_n)^mQ_nT(a_+^{-1})P_n
\quad(m\in\N).
\end{eqnarray*}
Taking into account (\ref{eq:BS1-6})--(\ref{eq:BS1-7}) and
Lemma~\ref{le:basic-Toeplitz}(b), (c), one can prove that the operators
$F_{n,m}$ can be written in the form (\ref{eq:BS1-3}).
\end{proof}
\begin{lemma}\label{le:BS2}
Under the assumptions of Lemma~{\rm\ref{le:BS1}},
\begin{eqnarray*}
&&
(P_n-P_{n-\kappa}) T_n^{-1}(a)P_{\kappa-1}
\\
&&=
(P_n-P_{n-\kappa}) T(a_-^{-1})(P_n-P_{n-\kappa}) T(b) P_{\kappa-1}T(a_-^{-1})P_{\kappa-1}
+X_{n,\kappa}
\end{eqnarray*}
and
\begin{eqnarray*}
&&
P_{\kappa-1}T_n^{-1}(a)(P_n-P_{n-\kappa})
\\
&&=
P_{\kappa-1}T(a_+^{-1})P_{\kappa-1}T(c)(P_n-P_{n-\kappa}) T(a_+^{-1})(P_n-P_{n-\kappa})
+Y_{n,\kappa},
\end{eqnarray*}
where
\begin{eqnarray}
\label{eq:defX}
X_{n,\kappa}&:=&
(P_n-P_{n-\kappa})H(a_+^{-1})H(\widetilde{c})Q_nT(b)P_{\kappa-1}T(a_-^{-1})P_{\kappa-1}
\\
\nonumber
&&
-(P_n-P_{n-\kappa})T(a_+^{-1})P_nT(c)
\\
\nonumber
&&
\quad
\times\sum_{m=1}^\infty\big(Q_nH(b)H(\widetilde{c})Q_n\big)^m
T(b)P_{\kappa-1}T(a_-^{-1})P_{\kappa-1},
\\
\label{eq:defY}
Y_{n,\kappa}&:=&
P_{\kappa-1}T(a_+^{-1})P_{\kappa-1}T(c)Q_nH(b)H(\widetilde{a_-^{-1}})(P_n-P_{n-\kappa})
\\
\nonumber
&&
-P_{\kappa-1}T(a_+^{-1})P_{\kappa-1}T(c)
\\
\nonumber
&&
\quad
\times
\sum_{m=1}^\infty\big(Q_nH(b)H(\widetilde{c})Q_n\big)^m
T(b)P_nT(a_-^{-1})(P_n-P_{n-\kappa}),
\end{eqnarray}
and the convergence in {\rm(\ref{eq:defX}), (\ref{eq:defY})} is understood in
the sense of $\cB(\ell^\Psi(\Z_+))$.
\end{lemma}
The first identity involving $X_{n,k}$ is actually proved in \cite[Section~8]{BS80}
in the setting of $\ell^2(\Z_+)$ (see also \cite[Theorem~3]{BS81} and
\cite[Theorem~2.2]{BG03}). An identity similar to the second
one was used in \cite[Theorem~4]{BS81}. A proof of the second identity is also
given  in \cite[Lemma~2.2]{K06c} in the setting of $\ell^2(\Z_+)$. Once we have
at hands Lemmas~\ref{le:basic-Toeplitz}, \ref{le:stability}, and \ref{le:BS1},
the proof of Lemma~\ref{le:BS2} can be developed by analogy with \cite[Section~8]{BS80}
or \cite[Lemma~2.2]{K06c} using the factorization technique of the proof of
Lemma~\ref{le:BS1}.
\subsection{The norms of $X_{n,\kappa}$ and $Y_{n,\kappa}$ are asymptotically small}
In this section we will show that the norms of of operators $X_{n,\kappa}$
and $Y_{n,\kappa}$ are asymptotically small whenever the generating function
$a$ is sufficiently smooth.

Put $\Delta_0:=P_0$ and $\Delta_j:=P_j-P_{j-1}$ for $j\in\{1,\dots,n\}$.
First we will prove the following auxiliary estimate for truncations of Toeplitz
operators.
\begin{lemma}\label{le:truncations}
Suppose $\Phi$, $\Psi$ are complementary $N$-functions both satisfying the
$\Delta_2^0$-condition and $\varphi=\{\varphi_k\}_{k=0}^\infty$,
$\psi=\{\psi_k\}_{k=0}^\infty$ are weight sequences in $\cW$.
There exists a constant $C>0$ depending only on $\Phi,\Psi$ and $\varphi,\psi$
such that if $a\in W\cap F\ell_{\varphi,\psi}^{\Phi,\Psi}$, then for every
$n\in\N$ and every $j\in\{0,\dots,n\}$,
\begin{eqnarray*}
\|Q_nT(a)\Delta_j\|_{\cB(\ell^\Psi(\Z_+))}
&\le&
C\frac{\|a-a^{(n-j)}\|_{F\ell_{\varphi,\psi}^{\Phi,\Psi}}}{\psi_{n-j+1}},
\\
\|\Delta_j T(a) Q_n\|_{\cB(\ell^\Psi(\Z_+))}
&\le&
C\frac{\|a-a^{(n-j)}\|_{F\ell_{\varphi,\psi}^{\Phi,\Psi}}}{\varphi_{n-j+1}}.
\end{eqnarray*}
\end{lemma}
\begin{proof}
This statement is proved similarly to \cite[Lemma~9]{KS05}. Suppose
$c=\{c_k\}_{k=0}^\infty\in\ell^\Psi(\Z_+)\setminus\{0\}$. Clearly,
\[
\Psi\left(\frac{|c_j|}{\|c\|_{\ell^\Psi(\Z_+)}}\right)
\le
\sum_{k=0}^\infty\Psi\left(\frac{|c_k|}{\|c\|_{\ell^\Psi(\Z_+)}}\right)
\le 1.
\]
Therefore,
\begin{equation}\label{eq:truncations-1}
|c_j|\le\Psi^{-1}(1)\|c\|_{\ell^\Psi(\Z_+)}.
\end{equation}
It is easy to check that
\begin{equation}\label{eq:truncations-2}
(Q_nT(a)\Delta_j c)_k=\left\{
\begin{array}{ll}
0, & 0\le k\le n,\\
a_{k-j}c_j, & k>n.
\end{array}\right.
\end{equation}
Without loss of generality assume that
$\|\{(a-a^{(n-j)})_k\}_{k\in\Z_+}\|_{\ell_\psi^\Psi(\Z_+)}>0$. Then taking into
account (\ref{eq:truncations-1}), (\ref{eq:truncations-2}), and that
$a_k=(a-a^{(n-j)})_k$ for $k>n-j$, we have
\begin{eqnarray*}
&&
\sum_{k=0}^\infty\Psi\left(
\frac{|(Q_nT(a)\Delta_jc)_k|\psi_{n-j+1}}
{\Psi^{-1}(1)\|a-a^{(n-j)}\|_{F\ell_{\varphi,\psi}^{\Phi,\Psi}}\|c\|_{\ell^\Psi(\Z_+)}}
\right)
\\
&&=
\sum_{k=n+1}^\infty\Psi\left(
\frac{|a_{k-j}|\psi_{n-j+1}}{\|a-a^{(n-j)}\|_{F\ell_{\varphi,\psi}^{\Phi,\Psi}}}
\cdot
\frac{|c_j|}{\Psi^{-1}(1)\|c\|_{\ell^\Psi(\Z_+)}}
\right)
\\
&&\le
\sum_{k=n+1}^\infty\Psi\left(
\frac{|a_{k-j}|\psi_{k-j}}{\|\{(a-a^{(n-j)})_k\}_{k\in\Z_+}\|_{\ell_\psi^\Psi(\Z_+)}}
\right)
\\
&&=
\sum_{k=n-j+1}^\infty\Psi\left(
\frac{|(a-a^{(n-j)})_k|\psi_k}{\|\{(a-a^{(n-j)})_k\}_{k\in\Z_+}\|_{\ell_\psi^\Psi(\Z_+)}}
\right)
\le 1.
\end{eqnarray*}
Therefore,
\begin{equation}\label{eq:truncations-3}
\|Q_nT(a)\Delta_j\|_{\cB(\ell^\Psi(\Z_+))}
\le
\frac{\Psi^{-1}(1)}{\psi_{n-j+1}}\|a-a^{(n-j)}\|_{F\ell_{\varphi,\psi}^{\Phi,\Psi}}.
\end{equation}
The second estimate is proved by using the duality argument. Obviously,
$Q_n^*=Q_n$, $\Delta_j^*=\Delta_j$, and $(T(a))^*=T(\overline{a})$.
Since $\Phi$ and $\Psi$ are complementary $N$-functions, by \cite[Section~13]{M83},
$[\ell^\Psi(\Z_+)]^*=\ell^\Phi(\Z_+)$ and
\begin{equation}\label{eq:truncations-4}
\|\Delta_jT(a)Q_n\|_{\cB(\ell^\Psi(\Z_+))}
\le
2\|Q_nT(\overline{a})\Delta_j\|_{\cB(\ell^\Phi(\Z_+))}.
\end{equation}
By Lemma~\ref{le:properties}(a), $\overline{a}\in F\ell_{\psi,\varphi}^{\Psi,\Phi}$ and
\begin{equation}\label{eq:truncations-5}
\|\overline{a}-\overline{a}^{(n-j)}\|_{F\ell_{\psi,\varphi}^{\Psi,\Phi}}
\le
C(\varphi,\psi,\Phi,\Psi)
\|a-a^{(n-j)}\|_{F\ell_{\varphi,\psi}^{\Phi,\Psi}}.
\end{equation}
Then applying (\ref{eq:truncations-3}) to $\overline{a}\in F\ell_{\psi,\varphi}^{\Psi,\Phi}$,
we get
\begin{equation}\label{eq:truncations-6}
\|Q_nT(\overline{a})\Delta_j\|_{\cB(\ell^\Phi(\Z_+))}
\le
\frac{\Phi^{-1}(1)}{\varphi_{n-j+1}}
\|\overline{a}-\overline{a}^{(n-j)}\|_{F\ell_{\psi,\varphi}^{\Psi,\Phi}}.
\end{equation}
Combining (\ref{eq:truncations-4})--(\ref{eq:truncations-6}), we arrive at
\[
\|\Delta_jT(a)Q_n\|_{\cB(\ell^\Psi(\Z_+))}
\le
\frac{2\Phi^{-1}(1)}{\varphi_{n-j+1}}C(\varphi,\psi,\Phi,\Psi)
\|a-a^{(n-j)}\|_{F\ell_{\varphi,\psi}^{\Phi,\Psi}}.
\]
The lemma is proved.
\end{proof}
The main result of this section is the following.
\begin{lemma}\label{le:small-XY}
Suppose $\Phi,\Psi$ are complementary $N$-functions both satisfying the
$\Delta_2^0$-condition and $\varphi=\{\varphi_k\}_{k=0}^\infty$,
$\psi=\{\psi_k\}_{k=0}^\infty$ are weight sequences in $\cW$. Let
\[
a_-^{\pm 1}\in(W\cap F\ell_{\varphi,\psi}^{\Phi,\Psi})\cap H_-^\infty,
\quad
a_+^{\pm 1}\in(W\cap F\ell_{\varphi,\psi}^{\Phi,\Psi})\cap H_+^\infty
\]
and put $a=a_-a_+$, $b=a_-a_+^{-1}$, and $c=a_-^{-1}a_+$. Then the norms
of the operators $X_{n,\kappa}$, $Y_{n,\kappa}$ defined by
{\rm (\ref{eq:defX}), (\ref{eq:defY})} satisfy
\[
\|X_{n,\kappa}\|_{\cB(\ell^\Psi(\Z_+))}=o(1/\psi_{n-\kappa}),
\quad
\|Y_{n,\kappa}\|_{\cB(\ell^\Psi(\Z_+))}=o(1/\varphi_{n-\kappa})
\quad (n\to\infty).
\]
\end{lemma}
\begin{proof}
By Lemma~\ref{le:algebra}(a), $a,b,c\in W\cap F\ell_{\varphi,\psi}^{\Phi,\Psi}\subset W$.
Hence, by Lemma~\ref{le:basic-Toeplitz}(a), the operator $H(b)H(\widetilde{c})$ is compact
on $\ell^\Psi(\Z_+)$. In view of Proposition~\ref{pr:strong}, the sequence of operators
$Q_n=I-P_n$ tends strongly to the zero operator on $\ell^\Psi(\Z_+)$, whence
$\|Q_nH(b)H(\widetilde{c})Q_n\|_{\cB(\ell^\Psi(\Z_+))}\to 0$ as
$n\to\infty$. Therefore, for all sufficiently large $n$,
\[
\left\|\sum_{m=1}^\infty(Q_nH(b)H(\widetilde{c})Q_n)^m\right\|_{\cB(\ell^\Psi(\Z_+))}
\le M_1(a_-,a_+)<\infty
\]
and
\begin{equation}
\label{eq:small-XY-1}
\|Y_{n,\kappa}\|_{\cB(\ell^\Psi(\Z_+))}
\le
M_2(a_-,a_+)\|P_{\kappa-1}T(c)Q_n\|_{\cB(\ell^\Psi(\Z_+))},
\end{equation}
where $M_1(a_-,a_+)$ is a positive constant depending only
$a_\pm$ and
\[
\begin{split}
M_2(a_-,a_+)
:= &
\|T(a_+^{-1})\|_{\cB(\ell^\Psi(\Z_+))}
\|H(b)H(\widetilde{a_-^{-1}})\|_{\cB(\ell^\Psi(\Z_+))}
\\
&-
\|T(a_+^{-1})\|_{\cB(\ell^\Psi(\Z_+))}
M_1(a_-,a_+)
\|T(b)\|_{\cB(\ell^\Psi(\Z_+))}
\|T(a_-^{-1})\|_{\cB(\ell^\Psi(\Z_+))}.
\end{split}
\]
On the other hand, by Lemma~\ref{le:truncations},
\begin{eqnarray}
\label{eq:small-XY-2}
\|P_{\kappa-1}T(c)Q_n\|_{\cB(\ell^\Psi(\Z_+))}
&\le&
\sum_{j=0}^{\kappa-1}\|\Delta_j T(c)Q_n\|_{\cB(\ell^\Psi(\Z_+))}
\\
\nonumber
&\le&
C\sum_{j=0}^{\kappa-1}
\frac{\|c-c^{(n-j)}\|_{F\ell_{\varphi,\psi}^{\Phi,\Psi}}}{\varphi_{n-j+1}}.
\end{eqnarray}
Since $\|c-c^{(n-j)}\|_{F\ell_{\varphi,\psi}^{\Phi,\Psi}}$ and
$1/\varphi_{n-j+1}$ are monotonically increasing with respect to $j$,
from (\ref{eq:small-XY-1}) and (\ref{eq:small-XY-2}) it follows that
\[
\|Y_{n,\kappa}\|_{\cB(\ell^\Psi(\Z_+))}
\le
M_2(a_-,a_+)C\kappa
\frac{\|c-c^{(n-\kappa)}\|_{F\ell_{\varphi,\psi}^{\Phi,\Psi}}}{\varphi_{n-\kappa}}.
\]
By Lemma~\ref{le:properties}(b),
$\|c-c^{(n-\kappa)}\|_{F\ell_{\varphi,\psi}^{\Phi,\Psi}}=o(1)$ as $n\to\infty$.
Thus $\|Y_{n,\kappa}\|_{\cB(\ell^\Psi(\Z_+))}=o(1/\varphi_{n-\kappa})$ as $n\to\infty$.
The equality $\|X_{n,\kappa}\|_{\cB(\ell^\Psi(\Z_+))}=o(1/\psi_{n-\kappa})$ as
$n\to\infty$ is proved analogously.
\end{proof}
Notice that one can get better estimates for the norms of $X_{n,\kappa}$ and
$Y_{n,\kappa}$ in the setting of the Hilbert space $\ell^2(\Z_+)$. In that
case effective estimates for the norms of $Q_nH(b)$ and $H(\widetilde{c})Q_n$
are easily available (see \cite[Section~10.35]{BS06} or \cite[Section~5]{K06a}).
Therefore one can conclude not only that
\[
\left\|\sum_{m=1}^\infty (Q_nH(b)H(\widetilde{c})Q_n)^m\right\|_{\cB(\ell^2(\Z_+))}=O(1)
\quad(n\to\infty)
\]
but that this norm tends to zero as $n\to\infty$ with some determined speed
depending on the smoothness of $a$. For instance, if $a\in W\cap F\ell_{\alpha,\beta}^{1,1}$
and $\alpha,\beta\ge 0$, $\alpha+\beta>1$, then $\|Q_nH(b)\|_{\cB(\ell^2(\Z_+))}=o(n^{-\beta})$,
$\|H(\widetilde{c})Q_n\|_{\cB(\ell^2(\Z_+))}=o(n^{-\alpha})$, and
\[
\left\|\sum_{m=1}^\infty (Q_nH(b)H(\widetilde{c})Q_n)^m\right\|_{\cB(\ell^2(\Z_+))}
=
o(n^{1-\alpha-\beta})
\quad(n\to\infty).
\]
This observation explains why the results of Theorem~\ref{th:main}
(based on estimates of of truncations of Toeplitz operators on $\ell^\Psi(\Z_+)$)
are less precise than the results of Theorems~\ref{th:nonvanishing-1} and
Theorems~\ref{th:nonvanishing-2} (based on estimates of
of truncations of Toeplitz and Hankel operators on $\ell^2(\Z_+)$).
\subsection{Proof of Theorem~\ref{th:main}(c)}
\label{subsec:proof}
We will prove the asymptotic formula for $D_n[a(t)t^\kappa]$ and $\kappa\in\N$
following \cite[Section~10.47]{BS06} (see also \cite[Section~3.3]{K06c}).
In view of Theorem~\ref{th:main}(a), there exist functions $a_-$ and $a_+$
such that $a=a_-a_+$ and
\[
a_-^{\pm 1}\in
(W\cap F\ell_{\varphi,\psi}^{\Phi,\Psi})\cap H_-^\infty
\subset
W\cap H_-^\infty,
\quad
a_+^{\pm 1}\in
(W\cap F\ell_{\varphi,\psi}^{\Phi,\Psi})\cap H_+^\infty
\subset
W\cap H_+^\infty.
\]
Thus all the conditions of Lemmas~\ref{le:stability}--\ref{le:BS2} and
\ref{le:small-XY} are fulfilled.
Therefore the matrices $T_n(a)$ are invertible for all sufficiently large $n$.
By Lemma~\ref{le:Jacobi},
\begin{equation}\label{eq:proof-1}
D_{n-\kappa}[a(t)t^\kappa]
=
(-1)^{n\kappa}D_n(a)\det[P_{\kappa-1}T_n^{-1}(a)(P_n-P_{n-\kappa})].
\end{equation}
By Lemmas~\ref{le:BS2} and \ref{le:small-XY},
\begin{eqnarray*}
&&
\|P_{\kappa-1}T_n^{-1}(a)(P_n-P_{n-\kappa})
\\
&&
-
P_{\kappa-1}T(a_+^{-1})P_{\kappa-1}T(c)(P_n-P_{n-\kappa})T(a_+^{-1})(P_n-P_{n-\kappa})
\|_{\cB(\ell^\Psi(\Z_+))}
\end{eqnarray*}
is $o(1/\varphi_{n-\kappa})$
as $n\to\infty$. Applying Proposition~\ref{pr:LA}
to the above matrices (which are of the size $m=\kappa$),
we get
\begin{eqnarray}
\label{eq:proof-2}
&&
\det[P_{\kappa-1}T_n^{-1}(a)(P_n-P_{n-\kappa})]
\\
\nonumber
&&=
\det[P_{\kappa-1}T(a_+^{-1})P_{\kappa-1}T(c)(P_n-P_{n-\kappa})T(a_+^{-1})(P_n-P_{n-\kappa})]
+o(1/\varphi_{n-\kappa})
\\
\nonumber
&&=
D_{\kappa-1}(a_+^{-1})\det[P_{\kappa-1}T(c)(P_n-P_{n-\kappa})]D_{\kappa-1}(a_+^{-1})
+o(1/\varphi_{n-\kappa})
\end{eqnarray}
as $n\to\infty$. From (\ref{eq:proof-1}) and (\ref{eq:proof-2}) we get
\begin{eqnarray}
\label{eq:proof-3}
D_n[a(t)t^\kappa]
&=&
(-1)^{(n+\kappa)\kappa}D_{n+\kappa}(a)D_{\kappa-1}^2(a_+^{-1})
\\
\nonumber
&&
\times\det[P_{\kappa-1}T(c)(P_{n+\kappa}-P_n)]+o(1/\varphi_n)
\quad
(n\to\infty).
\end{eqnarray}
By Theorem~\ref{th:KS},
\begin{equation}\label{eq:proof-4}
D_{n+\kappa}(a)=G(a)^{n+\kappa+1}E(a)\{1+o(1)\}
\quad(n\to\infty).
\end{equation}
Since $a_+^{-1}\in H_+^\infty$, the matrix $T_{\kappa-1}(a_+^{-1})$ is
triangular. It is not difficult to see that
\begin{equation}\label{eq:proof-5}
G(a)^\kappa D_{\kappa-1}^2(a_+^{-1})=\frac{G(a)^\kappa}{G(a_+)^{2\kappa}}=G(b)^\kappa.
\end{equation}
On the other hand,
\begin{equation}\label{eq:proof-6}
\det[P_{\kappa-1}T(c)(P_{n+\kappa}-P_n)]=\det\left(\begin{array}{ccc}
c_{-n-1} &\dots & c_{-n+\kappa-2} \\
\vdots & \ddots & \vdots\\
c_{-n-\kappa} & \dots &c_{-n-1}
\end{array}\right).
\end{equation}
Combining (\ref{eq:proof-3})--(\ref{eq:proof-6}), we arrive at the desired
formula
\begin{eqnarray*}
&&
D_n[a(t)t^{\kappa}]
=
G(a)^{n+1} E(a)\{1+o(1)\}
\\
&&
\quad\quad
\times
(-1)^{(n+\kappa)\kappa}G(b)^\kappa
\left\{\det\left(\begin{array}{ccc}
c_{-n-1} &\dots & c_{-n+\kappa-2}\\
\vdots & \ddots &\vdots\\
c_{-n-\kappa} &\dots & c_{-n-1}
\end{array}\right)+o(1/\varphi_n)
\right\}
\end{eqnarray*}
as $n\to\infty$. The proof of the asymptotic formula for $D_n[a(t)t^{-\kappa}]$
is analogous. \qed
\bibliographystyle{amsalpha}

\end{document}